\documentclass[12pt,a4paper]{article}

\usepackage{pdflscape}
\usepackage{amsmath}
\usepackage{amssymb}
\usepackage{latexsym}
\usepackage{graphics}
\usepackage{color}
\usepackage{epsfig}
\usepackage{color}
\usepackage{graphicx,epstopdf}

\newcommand{\re}{{\mathbb R}}

\newcommand{\cA}{{\mathcal{A}}}

\newcommand{\cT}{{\mathcal{T}}}
\newcommand{\cH}{{\mathcal{H}}}

\newcommand{\cP}{{\mathcal{P}}}
\newcommand{\cM}{{\mathcal{M}}}

\newcommand{\cJ}{{\mathcal{J}}}

\newcommand{\bx}{{\boldsymbol{x}}}
\newcommand{\by}{{\boldsymbol{y}}}

\newcommand{\bp}{{\boldsymbol{p}}}

\newcommand{\ba}{{\boldsymbol{a}}}
\newcommand{\bb}{{\boldsymbol{b}}}

\newtheorem{theorem}{Theorem}
\newtheorem{prop}{Proposition}
\newtheorem{lemma}{Lemma}
\newtheorem{cor}{Corollary}
\newtheorem{remark}{Remark}
\newtheorem{ex}{Example}
\newtheorem{defi}{Definition}

\date{}

\begin{document}

\author{Rinat Kamalov
\thanks{Gran Sasso Science Institute, Italy;
Moscow Institute of Science and Technology,  
{e-mail: \tt\small rinat020398god@yandex.ruDISIM}}\ ,  
Vladimir Yu. Protasov  
\thanks{University of L'Aquila, Italy
  {e-mail: \tt\small
vladimir.protasov@univaq.it}}}

\title{The length of switching intervals\\ of a stable linear system 
\thanks{Sections 1, 2, and 3 were written by V.Protasov; Sections 4 and 5 were written by 
R.Kamalov. All results
in this paper are products of authors collaborative work.
}
}

\maketitle

\begin{abstract}

The linear switching system is 
a system of ODE with the time-dependent  matrix  taking  values from a given 
control matrix set. The system is (asymptotically) stable if all its trajectories tend to zero
for every control function. We consider possible mode-dependent restrictions
on the lengths of switching intervals which keeps the stability of the system. 
When the stability of trajectories 
with short switching intervals implies the stability of all trajectories? 
To answer this question we introduce the concept of ``cut tail points'' of linear operators
and study them by the convex analysis tools. We reduce the problem to 
the construction of 
Chebyshev-type exponential polynomials, for which we derive an algorithm 
 and present the corresponding numerical results.

\bigskip

\smallskip

\begin{flushright}

\end{flushright}

\end{abstract}
\bigskip

\vspace{1cm}

\begin{center}

\large{\textbf{1. Introduction}}	
\end{center}
\bigskip 

Linear switching system is a linear ODE  $\dot \bx (t) \, = \,  A(t)\bx(t)$
on the vector-function $\bx \, : \, \re_+ \, \to \, \re^d$ with the initial 
condition~$\bx(0) = \bx_0$ such that  the matrix function $A(t)$ 
takes values from a given compact set~$\cA$ called 
 {\em control set}. The {\em control function}, 
 or the  
{\em switching law}  is an arbitrary measurable function  
$A: \re_+ \to \cA$. Each matrix~$A \in \cA$ is called a {\em regime} of the system.
 We shall identify the linear switching system with the corresponding 
control set~$\cA$. 
The {\em switching interval} is the segment between two 
consecutive switching points. We shall use the same word also for the length of this 
interval. 
 
Linear switching systems naturally appear  
in problems of robotics, electronic engineering, mechanics, planning, etc.~\cite{Liberzon}.
One of the main issues  is the stability of the system. 
A system is called {\em asymptotically stable} (in short,  ``stable''), if all its trajectories tend to zero. 
For one-element control sets~$\cA = \{A\}$, the stability is equivalent to that 
the matrix  $A$ is Hurwitz (or stable), i.e., all its 
eigenvalues have  negative real parts.  

 If $\cA$ contains more than one matrix, then the stability problem becomes 
much harder. In general, we  do not know the switching law 
that provides the maximal rate of growth of trajectories. Moreover, 
even some basic properties of those laws such that the frequency of switches or the 
lengths of switching intervals are unknown. In this paper we make one step towards the solution     
and prove that the ``worst'' trajectories do not have long switching intervals. 
More precisely, we show that for each matrix~$A \in \cA$, one can associate
a positive number~$T_{cut}(A)$, which is the moment of time 
when the trajectory of ODE generated by~$A$ enters the interior of its symmetrized 
convex hull. Our fist fundamental theorem states that 
if all the  trajectories whose switching intervals of the regime~$A$ do not exceed~$T_A$, then 
the system is stable. Therefore, if 
the trajectories for which the length of the 
switching intervals do not exceed the values~$T_{cut}(A)$ 
for the corresponding regimes~$A \in \cA$, then 
the system is stable.  Thus, the stability can be verified only for 
switching laws  with the restricted switching intervals. 
Thus, long switching intervals do not influence the stability. Moreover, we present a 
method of finding of the corresponding maximal lengths~$T_{cut}(A)$. We 
impose two assumptions: 
\smallskip 

\noindent \textbf{Assumption 1.} All the matrices~$A \in \cA$ are Hurwitz. 
\smallskip 

\noindent \textbf{Assumption 2.} There is a mode-dependent dwell time restriction: all switching intervals 
for every~$A \in \cA$ is at least~$m(A)$, where $m(A) >0$ is a given number.  
    
    \bigskip 
    
    Both those assumptions are natural. If one of the matrices form~$\cA$
    is not Hurwitz, then the system cannot be stable. Hence, it makes sense 
    to decide stability only under  Assumption~1. On the other hand, without 
    the dwell time assumption one can split all long switching intervals 
    to short parts by momentary interfering with other regimes. 
    In this case the length restrictions does not play a role, which 
    the problem becomes irrelevant.   
\smallskip 

Before formulating the main results we need to introduce one more concept,  
which is probably of some independent interest. 

\bigskip 

\begin{center}
\large{\textbf{2. Cut tail points}}
\end{center}
\bigskip

We are going to introduce the formal definition of the 
maximal switching interval~$T_A$ and observe some of its properties, after which 
we will be able to formulate the  theorem 
on the restriction of switching intervals for stable systems.

We consider a system~$\dot \bx\, = \, A\bx, \bx(0) = \bx_0$
with a constant Hurwitz  matrix~$A$,  without switches. We assume that~$\bx_0$  is a {\em generic point}, i.e., it does not belong to any proper invariant subspace of~$A$. The trajectory starting at~$\bx_0$ is also called generic. 
This does not reduce the generality because any trajectory can be made generic after 
reducing to the  corresponding invariant subspace. 

 The following fact is elementary, we include the proof for convenience of the reader. 
  \begin{lemma}\label{l.fill}
 A generic trajectory is not contained in a proper  
 subspace of~$\re^d$. 
 \end{lemma}  
 {\tt Proof.} Let $L$ be the linear span of the trajectory. Since 
  $\dot \bx(t) \in L$, we have~$A\bx(t) \in L$ for all~$t \ge 0$, and 
  therefore~$AL \subset L$. Thus, $L$ is an invariant subspace of~$A$ that contains 
  $\bx(0) = \bx_0$. Since~$\bx_0$ is a generic point, it follows that 
  $L = \re^d$.

{\hfill $\Box$}
\medskip 
 
 For a subset~$M \subset \re^d$,  we denote by ${\rm co}_s M\, =  \, {\rm co} \,\{M , -M\}$
 its {\em symmetrized convex hull}. 
  For a generic trajectory~$\bx(\cdot )$ and for a segment~$[t_1, t_2] \subset \re_+$, let~$\Gamma(t_1, t_2)\, = \,\bigl\{\bx(t) \ | \ 
t \in [t_1, t_2]\bigr\}$ and $G(t_1, t_2)\, = \, {\rm co}_s \, \Gamma(t_1, t_2)$.  
 It is seen easily that $G(t_1+h, t_2+h)\, = \, 
e^{hA}\, G(t_1, t_2)$. 
We also denote by  $\Gamma = \Gamma(0, +\infty)$ the entire trajectory and, respectively, 
$G = G(0, +\infty)$.  Thus, $G$ is the symmetrized convex hull of~$\Gamma$. 
The origin is an interior point of~$G$. Therefore, if the matrix~$A$ is 
Hurwitz, then~$G$ is a convex body, i.e., a compact set with a nonempty interior.   
 \begin{defi}\label{d.80}
 A point~$T>0$ is called a {\em cut tail point} for 
 the stable system~$\dot \bx\, = \, A\bx, \ \bx(0) = \bx_0$, if  
 for every $t> T$, the point $\bx(t)$  belongs to the interior of the set~$G(0, T)$. 
 \end{defi}   
 \begin{prop}\label{p.100}
For every stable system~$\dot \bx\, = \, A\bx, \ \bx(0) = \bx_0$, the set of cut tail points is nonempty. 
 \end{prop}  
 {\tt Proof.}  Since the point $\bx_0$ is generic, it follows that every arc of the trajectory is not contained in a hyperplane, therefore,  for each segment~$[t_1, t_2]$, 
the set $G(t_1, t_2)$
is full-dimensional and hence,  
it contains a ball centered at the origin. 
On the other hand, the trajectory converges to zero as~$t\to \infty$, hence, 
say, $G(0,1)$ contains $\bx(t)$ for all~$t$ bigger than some $T>1$. 
Since $G(0,1) \subset G(0,T)$, we see that~$T$ is a cut tail point.  

{\hfill $\Box$}
\medskip

The set  of cut tail points is obviously closed, hence, there exists
 the minimal point, which will be denoted by~$T_{cut} = T_{cut}(A)$. Clearly, 
$\bx(T_{cut})$ is on the boundary of~$G$. Thus, $T_{cut}$ is the  moment of time when the trajectory enters the interior of its 
symmetrized convex hull. To proceed further we need the following well-known auxiliary fact:  
if two trajectories 
start at generic  points~$\bx_0, \tilde \bx_0$, then 
they are affinely similar. To show this we apply the following 
\begin{lemma}\label{l.similarity}
For an arbitrary $d\times d$ matrix $A$ and for arbitrary points  $\ba, \bb \in \re^d$ 
such that  $\ba$ does not belong to an eigenspace of~$A$, there exists a $d\times d$ matrix~$C$ that 
commutes with  $A$ and takes  $\ba$ to $\bb$. 
\end{lemma}
{\tt Proof.}  
We pass to the Jordan basis of the matrix 
$A$ and consider one Jordan block~$\Lambda$ of size $s$ corresponding to an 
eigenvalue~$\lambda$. Denote by $\ba'$ and $\bb'$ the $s$-dimensional components of the 
vectors  $\ba$ and $\bb$ respectively corresponding to this block. Denote by $\cH$ the set of Hankel upper-triangular matrices of size~$s$, 
i.e., matrices of the form
$$
\left(
\begin{array}{ccccc}
\alpha_1 &  \alpha_2 & \cdots & \alpha_{s-1} & \alpha_{s} \\
0 &  \alpha_1 & \alpha_2 & \cdots & \alpha_{s-1} \\
0 &  0 & \alpha_1 & \alpha_2 & \cdots \\
\vdots & \vdots & \ddots & \ddots& \vdots \\
0 & 0  & \cdots & 0  & \alpha_1
\end{array}
\right)\, . 
$$
Note that the set $\cH$ is a linear space and is a multiplicative group. 
All elements of $\cH$ commute with $\Lambda$. The set $U \, = \, \{X\ba' \ | \ X \in \cH\}$ is a linear subspace of  
$\re^s$ invariant with respect to every matrix from~$\cH$. 
Let us prove that $U =\re^s$. If this is not the case, then~$U$ 
is an invariant subspace for each matrix from~$\cH$, in particular, for the matrix~$\Lambda$. 
 Therefore, $U$ coincides with one of spaces  
 $U_j = \bigl\{(x_1, \ldots , x_j, 0, \ldots , 0)^T \in \re^s\bigr\}$
 (they form the complete set of  invariant subspaces of~$\Lambda$). 
 Hence, $\bb' \in U_j$ for some~$j\le s-1$. Consequently, the vector $\ba$
 belongs to an eigenspace of the matrix~$A$ generated by all vectors of its Jordan basis 
 except for the last $s-j$ vectors, which correspond to the block~$\Lambda$. 
 This contradicts to the condition on~$\ba$. Thus, $U =\re^s$. 
 Hence, there exists a matrix $X \in \cH$ for which $X\ba = \bb$. 
  Let us remember that $X$ commutes with~$\Lambda$. 
  Having found such a matrix~$X$ for each Jordan block of the matrix~$A$
  we compose a block matrix~$C$ for which  $CA  = AC$ and $C\ba = \bb$.

{\hfill $\Box$}
\smallskip

Now we are ready  to establish the affine similarity of all 
generic trajectories. 
\begin{prop}\label{p.similarity}
All generic trajectories of a matrix~$A$ are affinely similar. 
\end{prop}
{\tt Proof.}  Let $\bx(\cdot), \by(\cdot)$ be two trajectories and 
their starting points~$\bx_0, \by_0$ do not belong to an invariant subspace of~$A$. 
Lemma~\ref{l.similarity} 
gives a matrix~$C$ such that  $CA = AC$ and $C\bx_0 = \by_0$. 
Then 
$$
C \dot \bx \ = \ C A \bx \ = \ A C \bx, \, \quad C\bx_0 \, = \, \by_0\, . 
$$
Thus, the function~$C\bx(\cdot)$ satisfied the same equation 
$\dot \by \, = \, A\by$ with the same starting point~$\by_0$, then by the uniqueness it coincides with~$\by(\cdot)$. We see that $\by(\cdot)\, = \, C\by(\cdot)$. 
It remains to note that  the similarity transform is nondegenerate. 
Otherwise the trajectory~$\by(\cdot)$ is not full-dimensional 
which contradicts to~Lemma~\ref{l.similarity}.

{\hfill $\Box$}
\smallskip

A direct consequence of Proposition~\ref{p.similarity} is that 
{\em the cut tail points do not depend on~$\bx_0$, provided~$\bx_0$ 
is generic}. Thus~$T_{cut}$ is  
a function of the matrix~$A$. In particular, speaking about 
the cut tail points we may not specify~$\bx_0$.

\begin{prop}\label{p.90}
For every generic trajectory, the arc~$\Gamma (0, T_{cut})$
lies on the boundary of~$G$.  
\end{prop}
{\tt Proof.} If, on the contrary, there exists~$t< T_{cut}$ such 
that~$\bx(t) \in {\rm int} \, G$, then, since the ``tail'' $\{\bx(s)\ | \ s> T_{cut}\}$
of the trajectory~$\Gamma$
lies in the interior of~$G$, 
it follows that~$\bx(t) \in {\rm int} \, G(0, T_{cut})$. 
Hence,
$$
\bx(T_{cut}) \ = \,  e^{(T_{cut} - t)A}\, \bx(t)\ 
\in \ {\rm int} \, e^{(T_{cut} - t)A}\, G(0, T_{cut})\ = \ {\rm int} \, 
G\, (T_{cut} - t\, , \, 2T_{cut}-t)\ \subset \ {\rm int}\, G.
$$
Thus, $\bx(T_{cut}) \in {\rm int}\, G$, which is a contradiction. 

{\hfill $\Box$}
\medskip

We see that the point  $\bx(T_{cut})$ separates the part of trajectory that lies on the 
boundary form that in the interior: $\bx(t) \in \partial G$, for all~$t\le T_{cut}$, and 
$\bx(t) \in {\rm int}\, G$, for all~$t > T_{cut}$. 
\begin{cor}\label{c.34}
The set of cut tail points is precisely the half-line
$[T_{cut}, +\infty)$. 
\end{cor}
\begin{cor}\label{c.36}
All extreme points of~$G$ are located on the union of the curves~$\Gamma(0, T_{cut})$
and $- \Gamma(0, T_{cut})$. 
\end{cor}

\bigskip 

\begin{center}
\large{\textbf{3. The first fundamental theorem}}
\end{center}
\bigskip

Now we are going to establish the relation 
between the cut tail points and the switching intervals for 
stable linear systems. Let us recall that we consider 
linear switching systems under Assumptions~1 and~2: 
all matrices~$A \in \cA$ are Hurwitz and the dwell time is bounded below. 
\begin{defi}\label{d.10}
Let~$\cA$ be a linear system   satisfying 
Assumptions~1 and~2 and let~$\cJ$ be an arbitrary subset of~$\cA$. 
The system is called~$\cJ$-stable if 
all its trajectories tend to zero provided for every~$A \in \cJ$  all switching
intervals corresponding to~$A$ do not exceed~$m(A) + T_{cut}(A)$. 
\end{defi}
Thus, a $ \cJ$-stable system is stable not for all trajectories but only for those 
with restricted lengths of intervals for the  regimes from~$\cJ$. 
In case $\cJ = \emptyset$ the~$\cJ$ stability becomes the usual stability. 
 In the proof we use the following result 
 from~\cite{PK}. To every~$A \in \cA$ we associate 
 numbers~$m(A) > 0$ and $M(A) \in (m(A), +\infty]$. 
 Denote by $\cM$ the set of pairs $\bigl(m(A), M(A)\bigr)\, , \ A \in \cA$. 
 We call a trajectory~{\em admissible} if for every~$A \in \cA$,  the lengths of all switching
intervals corresponding to~$A$ lie on the segment~$\bigl[m(A), M(A) \bigr]$.
 The system is called~$\cM$-stable if all its admissible trajectories 
 tend to zero.  To every 
$A \in \cA$, we associate a norm $f_A$ on~$\re^d$. 
The collection~$f=\{f_A\ | \ A \in \cA\}$ is called  a~{\em multinorm}. 
We assume that those norms are uniformly bounded. 
The following theorem was proved in~\cite{PK}. 
\smallskip 

\noindent \textbf{Theorem A}. {\em A system is~$\cM$-stable if and only if 
there exists a multinorm~$f$ that possesses the 
following property: for every admissible  trajectory~$\bx(\cdot)$ and for an arbitrary 
its switching point~$\tau$ between some regimes~$A$ and~$A'$, 
we have~$f_A(\bx(\tau))\, > \, \sup_{t \in [\tau + m(A'), \tau + T]} f_{A'}(\bx(t))$, 
where $m=m(A)$ and $\tau + T$ is the next switching point (or~$T=+\infty$ if $\tau$ is the largest switching point)}. 
\smallskip 

Now we formulate our first fundamental theorem on switching interval restrictions. 
\begin{theorem}\label{th.10}
Let a linear switching system~$\cA$ satisfy Assumptions~1 and~2.
Then for every subset~$\cJ \subset \cA$,  the~$\cJ$-stability
of~$\cA$ implies 
the stability.  
\end{theorem}
{\tt Proof}. For~$A \in \cJ$,  we set~$M(A) = m(A) + T_{cut}(A)$, and 
$M(A) =  +\infty$ for all~$A\notin \cJ$. If the system is~$\cJ$-stable, then 
it is~$\cM$-stable, and by Theorem~A there exists a multinorm~$f$
such that~$f_A(\bx(\tau))\, > \, \sup_{t \in [\tau + m(A'), \tau + T]} f_{A'}(\bx(t))$
for every switching point~$\tau$ from a regime~$A$ to~$A'$, 
where  $T \le M(A')$ is the length of the switching interval.  
Note that since the function~$f_{A'}$ is convex, it follows that the 
its supremum over the arc~$\gamma(\tau + m(A'),\tau + T )\, = \, \{\bx(t) \ | \ \tau + m(A') \, \le \, 
 \tau + T \}$ is equal to the supremum over the symmetrized convex 
 hull~$G(\tau + m(A'),\tau + T )$ of that arc. 
On the other hand, by the definition of~$T_{cut}$, for all~$T > T_{cut}(A)$, 
we have~$G(\tau + m(A'),\tau + T )\, = \, G(\tau + m(A'),\tau + T_{cut} )$. 
Hence, if the system is~$\cM$-stable for~$M(A) = m(A) + T_{cut}(A)$, 
then the assumptions of Theorem~A are satisfied for
all larger $M(A)$ and therefore, for~$M(A) = +\infty$. Now by the same~Theorem~A
applied for~$M(A)= +\infty, A\in \cA$,  we conclude that 
the system is stable.   
 
{\hfill $\Box$}
\smallskip

The multinorm from Theorem~A is called the {\em Lyapunov multinorm}
for~$\cM$-stability. 
Thus, if the system with the lengths of all switching intervals 
bounded by~$T_{cut}(A)$ is stable, then it is stable 
without any restrictions. In the proof of Theorem~\ref{th.10} we established that 
not only the stability but also the Lyapunov multifunction 
stays the same after removing all the restrictions. 

To restrict the lengths switching intervals in applied problems we need to 
find or estimate the cut tail point for given matrices. 
The definition of cut tail points involves the convex hull of the trajectory. 
Computation of a convex hull in~$\re^d$  is notoriously hard, therefore,~$T_{cut}$
 can unlikely be found in a straightforward way. 

\bigskip 

\begin{center}
\large{\textbf{4. Computing the cut tail points and the switching intervals}}
\end{center}
\bigskip

 To find~$T_{cut}$ for a given Hurwitz matrix~$A$ one needs to decide 
 if  a given~$T$ is a cut tail point?  We reduce this problem to finding the 
 quasipolynomial of be,            st approximation and present an algorithm for that. 
 The numerical results  are reported in the next section.  
 
 Let us denote by 
$\cP_{A}$ the linear span of the functions $f_{\bx}(t) = e^{tA}\bx, \  \bx \in \re^d$. 
This is the  space of quasipolynomials 
which are linear combinations of functions~$t^ke^{\, \alpha t}
\cos \beta t$ and $t^ke^{\, \alpha t}\, \sin \beta t\, $, 
where $\alpha + i \beta$ is an eigenvalue of~$A$ and $k = 1, \ldots , r-1$, where $r$ is the size of the 
largest Jordan bock of~$A$ corresponding to that eigenvalue.  The dimension of~$\cP_A$
is equal to the degree of the minimal annihilating polynomial of the matrix~$A$. 
\begin{theorem}\label{th.70}
Let~$A$ be a Hurwitz matrix and~$T>0$ be a number; then $T > T_{cut}$ if and only if  
the value of the 
convex extremal problem 
\begin{equation}\label{eq.extr-cut}
\left\{
\begin{array}{l}
p(T) \to \max \\
\|\bp\|_{C(\re_+)} \ \le \ 1\\
p \in \cP_A
\end{array}
\right. 
\end{equation}
is smaller than~$1$.  
\end{theorem}
{\tt Proof}. By Proposition~\ref{p.90}, the inequality~$T> T_{cut}$ is equivalent to that 
$\bx(T) \in {\rm int}\, G$. This means that the point $\bx(T)$
cannot be separated from~$G$ by a linear functional, i.e.,  
for every nonzero vector~$\bp\in \re^d$, we have 
$\bigl(\bp , \bx(T)\bigr) \, <   \, \sup\limits_{\by  \in G} \, (\bp , \by)$. 
Since~$G$ is a convex hull of points $\, \pm \, \bx(t)\, , \, t \ge 0$, 
it follows that $\sup\limits_{\by  \in G} \, (\bp , \by)\, = \, 
\sup\limits_{t \ge 0} \, \bigl|(\bp , \bx(t))\bigr|$. Thus,  
$$
\bigl(\bp , \bx(T)\bigr) \ <   \ \sup\limits_{t \in \re_+}\bigl| (\bp , \bx(t)) \bigr|.
$$ 
On the other hand, $(\bp , \bx(t)) = p(t)$, where $p \in \cP_{A}$
is the quasipolynomial with the vector of coefficients~$\bp$.  
We obtain $p(T) < \|p\|_{C(\re_+)}$. 
 Consequently, $p(T)< 1$
for every quasipolynomial from the unit ball~$\|p\|_{C(\re^+)} \le 1$.
Now by the compactness argument we conclude that the value of the 
problem~(\ref{eq.extr-cut}) is smaller than one.

{\hfill $\Box$}
\medskip

\begin{theorem}\label{th.80}
Suppose~$A$ is a stable matrix; then a number~$T$ is a cut tail point 
for~$A$ if and only if there exists a quasipolynomal~$p \in \cP_{A}$
for which~$T$ is a point of absolute maximum on~$\re^d$. 
\end{theorem}
{\tt Proof}. Proposition~\ref{p.90} implies that the point~$\bx(T)$ lies on the 
boundary of~$G$ precisely when~$T\le T_{cut}$. 
On the other hand, this is equivalent to say that~$\bx(T)$ can be separated from~$G$
by some linear functional. Arguing as in the proof of Theorem~\ref{th.70}
we obtain a quasyipolynomial~$p$ such that~$p(T) \ge  \|p\|_{C(\re_+)}$, i.e., 
$T$ is a point of maximum of~$p$ on~$\re^+$. Thus,~$T\le T_{cut}$
if and only if there exists a quasipolynomal~$p$ that has maximum at the point~$T$.

{\hfill $\Box$}
\medskip 

The problem~(\ref{eq.extr-cut}) is convex  and can be solved by the 
 convex optimisation  methods.  However, its main disadvantage that it contains 
 the norm of a quasipolynomial on the whole real line, which can cause numerical problems. 
 Theorem~\ref{th.90} below  is the second main result of this work. 
 It reduces problem~(\ref{eq.extr-cut}) to the following  problem on a compact interval: 
 \medskip 

\begin{equation}\label{eq.extr-cut1}
\left\{
\begin{array}{l}
\|\bp\|_{[0,T]} \ \to \ \min\\
p(T) = 1 \\
p \in \cP_A
\end{array}
\right. 
\end{equation}
Clearly, the answer to this problem cannot be smaller than one. 
Theorem~\ref{th.90} claims that it is equal to one precisely 
when~$T \ge T_{cut}$. 

\begin{theorem}\label{th.90}
A number~$T>0$ is a cut tail point 
if and only the value of the problem~(\ref{eq.extr-cut1}) is equal to one. 
\end{theorem}
{\tt Proof}. If for some~$p \in \cP_{\cA}$, we have~$\|\bp\|_{[0,T]} = 1$, then 
the largest point of maximum of~$p$ is bigger than or equal to~$T$, which in view of Theorem~\ref{th.80} means that~$T_{cut} \ge T$. Conversely, if the answer 
to the problem~(\ref{eq.extr-cut1}) is larger than one, 
then none of~$p \in \cP_{\cA}$ can have maximum in~$T$ and hence~$T > T_{cut}$.

{\hfill $\Box$}
\medskip

 If we are able to solve the problem~(\ref{eq.extr-cut1}) for every~$T$, then $T_{cut}$
 can be found merely by double division. 
\bigskip 

Thus, for given~$T>0$, we need to solve problem~(\ref{eq.extr-cut1}) of minimizing 
the norm of the quasipolynomial~$p \in \cP_A$ on the segment~$[0,T]$
under the constraint~$p(T) = 1$. Since such quasipolynomials form 
an affine  space of dimension~$d-1$ and the function~$f(p) = \|p\|_{C[0,T]}$
is convex, it follows that we can invoke the refinement theorem~(see, for instance,~\cite{MIT}), 
and conclude that there exist at most $d$ points~$t_1 < \ldots < t_n$, where~$\ n\le d$ on the segment~$[0,T]$
such that 
$$
\min_{p\in \cP_A, \, p(T) = 1}\|p\|_{C[0, T]}\ = \ \min_{p\in \cP_A, \, p(T) = 1}\, \max_{i=1, \ldots , n}|p(t_i)|\, . 
$$ 
For the optimal quasipolynomial~$\bar p$, which exists by the coercivity of the function~$f$, we must have~$|\bar p(t_i)| = \|\bar p\|, \, i = 1, \ldots , n$. For all other feasible quasipolynomials~$p$, 
the absolute value of~$p(t_i)$ is larger than~$\|\bar p\|$ for at least one point~$t_i, 
\, i=1, \ldots, n$. 
To solve the problem~(\ref{eq.extr-cut1}) and to find the best quasipolynomial~$\bar p$
we suggest the following algorithm. 
\medskip 

\noindent \textbf{Algorithm~1}. {\em Initialisation.} We have the space of  quasipolynomials 
$\cP_A$ of degree~$d$ on the segment~$[0,T]$. 
Take arbitrary $d$ points on the segment~$[0,T]$, denote this set by~$\cT_0$
and choose arbitrary admissible~$p_0\in \cP_A$ (i.e., $p_0(T) = 1$). 
Set~$b_0=1, B_0 = \|p_0\|$. Choose arbitrary~$\varepsilon > 0$. 
\smallskip 

{\em Main loop.} The $k$th iteration. We have a set of~$n \le d$ points~$\cT_{k-1}$
and the numbers $b_{k-1}, B_{k-1}$ which are 
the lower and upper bounds respectively for the value of the problem~(\ref{eq.extr-cut1}). 
Solve the linear programming problem
 \begin{equation}\label{eq.extr-cut2}
\left\{
\begin{array}{l}
r \ \to \ \min\\
|p(t_j)| \le r, \ t_j \in \cT_{d-1} \\
p \in \cP_A, \ p(T) = 1. 
\end{array}
\right. 
\end{equation}
The value~$r$ of this problem is a lower bound for the value of the 
problem~(\ref{eq.extr-cut1}).  Hence,~$b_{k} = \max\{b_{k-1}, r\}$
is the new lower bound. If~$\bar p$ is an optimal quasipolynomial for the 
problem~(\ref{eq.extr-cut2}), then we find the point~$s$
of absolute maximum of $|\bar p(t)|$ on the segment~$[0,T]$.  
Then $|\bar p(s)|$ is an upper bound for the value of the 
problem~(\ref{eq.extr-cut1}).  Hence,~$B_{k} = \min\{B_{k-1}, |\bar p(s)|\}$
is the new upper bound. 

We add the point~$s$ to the set~$\cT_{d-1}$ and remove those points~$t_i$ 
for which~$|\bar p(t_i)| < r$. 
 Thus, we obtain a new set~$\cT_{k}$. 
Then pass to the next iteration. 
 \smallskip 

{\em Termination.} If~$B_k - b_k < \varepsilon$, then stop. 
The value of the problem~(\ref{eq.extr-cut1}) belongs to the segment~$[b_k, B_k]$
\bigskip 
 
 \begin{remark}\label{r.25}
 {\em If all exponents of the quasipolynomials are real, i.e., 
 if the spectrum of the matrix~$A$ is real, 
 then these exponents form a Chebyshev system~\cite{Dz, KS}.  Hence, the optimal polynomial 
 in the problem~(\ref{eq.extr-cut1}) has an alternans at the points~$\cT_k$, 
 which can be efficiently found by the Remez algorithm~\cite{Rem}. In the general case, however, 
 the system is not Chebyshev and the Remez algorithm is not applicable.
 Algorithm~1 can be considered as a kind of replacement 
 of the Remez algorithm for non-Chebyshev systems.} 
 \end{remark}
 Numerical results presented in Section~5 show that 
 Algorithm~1 can effectively find decide whether~$T < T_{cut}$ and hence, by the double division, find~$T_{cut}$.  
 \bigskip 
 
 For $2\times 2$ matrices, the parameter $T_{cut}$ can be evaluated geometrically and in an  explicit form. 
 \bigskip 
 
 \noindent \textbf{The two-dimensional systems}. 

\smallskip 

We have a real $2\times 2$-matrix~$A$. 
For the sake of simplicity we exclude the case of multiple eigenvalues. 
So, we assume that~$A$  either has two different real eigenvalue or 
two complex conjugate eigenvalues. 
\smallskip  

\noindent {\tt Case 1. The eigenvalues of $A$ are real and different}. 
Denote them by $\alpha_1, \alpha_2$. 
 The trajectory~$\Gamma =  \{\bx(t)\, | \, t \ge 0\}$ in the basis of eigenvectors of~$A$ 
has the equation~$(x_1(t),  x_2(t)) \, = \, 
\bigl( e^{\alpha_1 t} \, , \,   e^{\alpha_2 t}\bigr)\, , \ 
t \in \re_+$. 
Since $A$ is stable, both~$\alpha_1, \alpha_2$ are negative 
and~$\bx(0) = (1, 1), \, \bx(+\infty) = (0, 0)$. 
Let~$\ba \in \Gamma$ be the point where the tangent line 
drawn from the point~$(-1, -1)$ touches~$\Gamma$. 
Then $G \, = \, {\rm co}_s \, \Gamma$ is bounded by the line segment 
from $(-1, -1)$ to $\ba$ and by the arc of~$\Gamma$ from $\ba$ to 
$(1, 1)$, then reflected about the origin. 
Hence, if $\ba = \bx(T)$, then $\bx(t) \in {\rm int}\, G$ for all~$t > T$, and 
so $T = T_{cut}$. We have~$\, \ba \, + \, s \, \dot \bx(T_{cut}) \, = 
\, (-1,-1)$, where $s> 0$ is some number. 
Writing this equality coordinatewise, we obtain the equation for~$t = T_{cut}$: 
\begin{equation}\label{eq.prec1}
\left\{
\begin{array}{lcl}
e^{\alpha_1 t}\ + \ s\, \alpha_1\, e^{\alpha_1 t}& = & - 1\\
e^{\alpha_2 t}\ + \ s\, \alpha_2\, e^{\alpha_2 t}& = & - 1
\end{array}
\right. 
\end{equation}
which becomes after simplification 
$\frac{1 + e^{- \alpha_1 t}}{\alpha_1}\, = \, \frac{1 + e^{- \alpha_2 t}}{\alpha_2}$. 
The unique solution is~$T_{cut}$. 
\medskip

{\tt Case~2. The eigenvalues of $A$ are complex conjugate}. They can be written as 
~$\alpha \pm i \beta$ with $\alpha < 0$. 
 The trajectory~$\Gamma$ in a suitable basis 
 has the equation~$(x_1(t), x_2(t)) \,  = \, 
 e^{\alpha t}\, (\cos \beta t \, , \, \sin \beta t)
 \, , \ 
t \in \re_+$. The trajectory $\Gamma$
 goes from the point~$\bx(0) = (1, 0)$ to zero making infinitely  many rotations. 
Taking the point of tangency~$\ba$ of $\Gamma$ with the line 
going from the point~$(-1, 0)$ and arguing as above we obtain 
\begin{equation}\label{eq.prec2}
\left\{
\begin{array}{lcr}
e^{\alpha t} \cos \beta t \ + \ s\, e^{\alpha t}\bigl(\alpha \cos \beta t \, - \, 
\beta \sin \beta t \bigr)& = & - 1\\
e^{\alpha t} \sin  \beta t \ + \ s\, e^{\alpha t}\bigl(\alpha \sin \beta t \, + \, 
\beta \cos \beta t \bigr)& = & 0
\end{array}
\right. 
\end{equation}
from which it follows  $\alpha \sin \beta t \, + \, \beta \cos \beta t\, + 
\, \beta \, e^{\alpha t} \, = \, 0$. The unique solution of this equation is $t = T_{cut}$. 
 
 \newpage 
 
 \bigskip 

\begin{center}
\large{\textbf{5. Numerical results}}
\end{center}
\bigskip
 
 In this section we present results of computation of $T_{cut}$
 in several examples of dimensions  from~$2$ to $4$. In dimension $2$ 
 we did the computation in two ways~by solving functional equations from the last section 
 and by Algoritj~1. We will see that the results are very close.
 In higher dimensions  we applied Algorithm~1. All computations take a few seconds in a standard laptop. 
 
\begin{figure}[h!]
 	{\includegraphics[scale = 0.6]{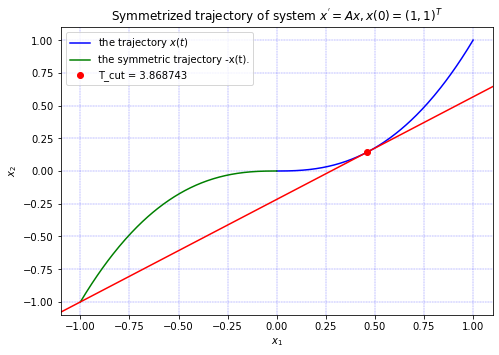}}
 	\caption{}
 	\label{picture:real}
 \end{figure}

 \begin{ex}
{\em 	 For the matrix $A =  \left( \begin{array}{cc}
	 		-0.2 & 0 \\
	 		0 & -0.5
	 	\end{array} \right)$, 
	we have  $ \lambda_1 = -0.2; \lambda_2 = -0.5 $. 
	The space of quasipolynomials is $ \mathcal{P}_A = \textrm{Lin}\{ e^{-0.2 t}, e^{-0.5 t}\} $. Algorithm~1 produces the value $ T_{cut} = 3.873191 $. 
	The value of~$T_{cut}$ found by solving the equation~(\ref{eq.prec1}) is $3.868743$, see fig. $ \ref{picture:real} $). }

\end{ex}

\begin{figure}[h!]
	{\includegraphics[scale = 0.6]{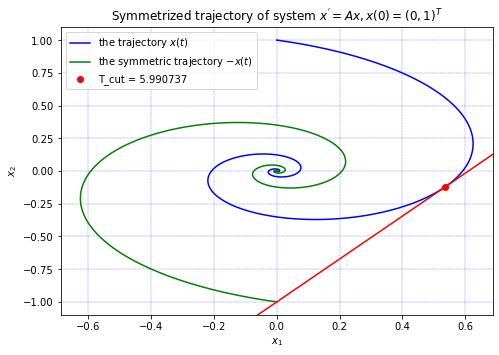}}
	\caption{}
	\label{picture:complex}
\end{figure}

\begin{ex}
	{\em  For the matrix $A =  \left( \begin{array}{cc}
	 		-0.1 & -0.3 \\
		 0.3 & -0.1
	 	\end{array} \right)$, 
	we have  $ \lambda_{1,2} \, = \, \alpha \pm \beta i = -0.1 \pm 0.3 i  $. 
	The space of quasipolynomials is $ \mathcal{P}_A = \textrm{Lin}\{ e^{-0.1 t}  \sin(0.3 t), e^{-0.1 t} \cos(0.3 t) \} $. Algorithm~1 produces the value $ T_{cut} = 5.9952163  $. 
	The value of~$T_{cut}$ found by solving the equation~(\ref{eq.prec1}) is $5.990737 $, see fig. $ \ref{picture:complex} $}

\end{ex}

\begin{ex}
	{\em For the  ${4 \times 4} $ matrix $ A \, = \, \left( \begin{array}{cccc}
		4.8 & 7.5 & 16 & 12 \\
		 -3 & -4.7 & -9.6 & -7.2 \\
		 -3.2 & -4.8 & -5.4 & -3.6 \\
		 4 & 6 & 6 & 3.9 \\
	\end{array} \right), $ 
	we have  $ \lambda_1 = -0.1; \lambda_2 = -0.2; \lambda_2 = -0.5; \lambda_2 = -0.6 $. The 
	corresponding space of quasipolynomials is $ \mathcal{P}_A = \textrm{Lin}\{ e^{-0.1 t}, e^{-0.2 t}, e^{-0.5 t}, e^{-0.6 t} \} $.  Algotithm~1 finds 
	 $ T_{cut} =  17.75795 $. Thus, if the matrix~$A$ presents in the 
	 linear switching systems, then the stability analysis can be done only 
	 for switching laws with the intervals of the regime~$A$ shorter than~$\, 17.75795 $. }
\end{ex}	

\begin{ex}
	{\em For the ${4 \times 4} $ matrix  $A \, =  \, \left( \begin{array}{cccc}
		22.8 & 15.7 & -74.5 & -41.9 \\
		-11.6 & -8.6 & 35.8 & 20.4 \\
		-1.2 & -0.3 & 5.8 & 3 \\
		 10.2 & 5.7 & -41.4 & -21.2 \\
	\end{array} \right), $ the spectrum is  $ \alpha_1 \pm \beta_1 i = -0.1 \pm 0.7; \alpha_2 \pm \beta_2 i = -0.5 \pm 0.3 i $. The corresponding space is $ \mathcal{P}_A = \textrm{Lin}\{ e^{-0.1 t}  \sin(0.7 t), e^{-0.1 t} \cos(0.7 t), e^{-0.5 t}  \sin(0.3 t) $, \\ $ e^{-0.5 t} \cos(0.3 t) \} $.  Applying algorithm~1 we get $ T_{cut} =  8.94363 $.
	Thus, for linear switching systems that include this matrix, 
	all the corresponding switching intervals can be chosen smaller than $8.94363 $}
\end{ex}	

\begin{ex}
{\em 	For the  matrix $\left( \begin{array}{cccc}
		-1.24 & -0.29 & -0.58 & -0.69 \\
		 -0.765 & -1.515 & -1.105 & -0.415 \\
		  1.31 & 1.66 &  0.37  & 0.86 \\
		  0.635 & -0.315 & 0.695 & 0.185  \\
	\end{array} \right),$ the spectrum is  $ \lambda_{1,2} = -0.3; \alpha_1 \pm \beta_1 i = -0.8 \pm 0.9 i $. We have  $ T_{cut} =  7.09526 $.}
\end{ex}		
\vspace{1cm}  

\noindent \textbf{Acknowledgements}. The work of Vladimir Yu.Protasov was supported by the Russian Science Foundation under
grant 20-11-20169 and was performed in Steklov Mathematical Institute of Russian Academy of Sciences.

 \end{document}